\documentclass[runningheads]{llncs}
\usepackage[T1]{fontenc}
\usepackage{graphicx}
\begin{document}
\title{Numerical modeling of elastic waves in thin shells with grid-characteristic method}

\author{Beklemysheva K.A.\inst{1}\orcidID{0000-0003-2769-6688} \and
Mikhel E.A.\inst{1}\orcidID{0009-0009-7224-9458} \and
Ovsyannikov A.D.\inst{1}\orcidID{0009-0008-3403-5611}}

\authorrunning{Beklemysheva K. et al.}
\titlerunning{Numerical modeling of elastic waves in thin shells with GCM}

\institute{
Moscow Institute of Physics and Technology, 9 Institutskiy per., Dolgoprudny, Moscow Region, 141701, Russian Federation}

\maketitle 
\begin{abstract}
Numerical modeling of strength and non-destructive testing of complex structures such as buildings, space rockets or oil reservoirs often involves calculations on extremely large grids. 
The modeling of elastic wave processes in solids places limitations on the grid element size because resolving different elastic waves requires at least several grid elements for the characteristic size of the modeled object. 
For a thin plate, the defining size is its thickness, and a complex structure that contains large-scale thin objects requires a large-scale grid to preserve its uniformity. 
One way to bypass this problem is the theory of thin plates and shells that replaces a simple material model on a fine three-dimensional mesh with a more complex material model on a coarser mesh. 
This approach loses certain fine effects inside the thin plate, but allows us to model large and complex thin objects with a reasonable size calculation grid and resolve all the significant wave types. 
In this research, we take the Kirchhoff-Love material model and derive a hyperbolic dynamic system of equations that allows for a physical interpretation of eigenvalues and eigenvectors. 
The system is solved numerically with a grid-characteristic method. Numerical results for several model statements are compared with three-dimensional calculations based on grid-characteristic method for a three dimensional elasticity.

\end{abstract}

\section{Introduction}

Elastic waves in thin plates and shells exhibit complex behavior and require specific approaches that are still under development. 
Generally, elastic waves can be categorized into volumetric waves (longitudinal and transverse) and surface waves \cite{aki_richards}. 
The most widespread applications like sonars and medical ultrasound are based on propagation and reflection of volumetric waves. 
Surface waves, even if they emerge in a specific problem statement, are considered weak and filtered as noise. 
In other applications of elastic waves diagnostics, like seismology or certain types of technical ultrasound, surface waves are considered separately and specifically.
 
Classic approaches to the study of surface waves and terms are generalized in \cite{surface_waves}. 
The most interesting types for the technical ultrasound are the waves that appear on a free surface of a solid body, Rayleigh waves and Lamb waves. 
Rayleigh waves appear in an elastic half-space and Lamb waves appear in a flat plate with two free surfaces. 
If the thickness of the plate increases to the infinity, Lamb waves turn into Rayleigh waves. 
Current studies show that Lamb waves are still of interest both to mathematicians and to actual engineering problems \cite{lamb_nonlinear}. 
The numerical modeling in this research is devoted to nonlinear Lamb waves propagating along a discontinuous plate and calculations are performed with the transient finite element method. 
Analytical approaches in this area are rare, but they still exist even for complex statements with multilayered plates \cite{lamb_analytical}. 
The results can be applied to metal and metal-ceramic composites that use only isotropic layers. 
A numerical end experimental study of negative refraction and focusing of surface waves is described in \cite{lamb_experimental}, where a step change in thickness is considered. 
This concise problem statement is useful both theoretically and practically.

The study on Lamb waves in a thin shell submerged in an acoustic medium is described in \cite{lamb_in_water}. 
It focuses on the backscattering problem and the analysis of the sound wave pattern reflected from an air-filled thin metal shell. 
Results are obtained in form of dispersion curves, modes and resonances. 
This type of waves is called leaky Lamb waves -- they propagate in thin shells that are in contact with an acoustic medium that is dense enough for acoustic energy to leak into it in a significant amount.
In \cite{lamb_in_air} a statement with a thin aluminum layer with air on both sides is considered. 
Experimental and numerical results are in a good compliance with each other and show dispersion and attenuation of leaky Lamb waves. 
\cite{lamb_in_21cent} describes a specific type of leaky Lamb waves and their influence on sound scattering patterns.
\cite{lamb_MSU} is concentrating on energy properties of surface waves in a submerged plate.

Another engineering problem that involves thin plates is the barely visible impact damage (BVID) in aircraft composites \cite{abrate_bvid}. The problem of BVID (damage propagation during impact, residual strength and non-destructive testing) is studied with a wide variety of approaches. A review of experimental studies is given in \cite{clifton_bvid} and a review of numerical models is given in \cite{bvid_numerical}.

The shell theory allows to significantly accelerate the calculations. It is widely used in the analysis of the postbuckling composites behavior and residual strength analysis \cite{bvid_residual}. 
In \cite{shells_vs_czm} two numerical models are compared, cohesive zone method (CZM) and Solid-Shell elements. CZM is more computationally expensive because it requires to model all composite plies separately while the Solid-Shell method generalizes all plies to a complex homogeneous medium. Results of the study show that the Solid-Shell method does not capture the delamination itself as well as the CZM, but shows similar results in terms of predicting consequences of complex loads. 
Shell material models are implemented in commercial software. In \cite{shells_fem_delamination} the effects of delaminations on the buckling loads of composite shells are studied with four-node linear shell elements with reduced integrals (S4R) implemented in ABAQUS.
An inverse shell finite element method is described in \cite{shells_inverse} and it involves a more complex material model. It is based on the refined zigzag theory \cite{shells_zigzag} that represents a more complex approach to through-thickness behavior. In \cite{shells_inverse} a time-domain modeling is performed for cantilevered composite plates under different harmonic loads. The distribution of the membrane damage index based on the equivalent von Mises strain is analyzed for flat and curved plates.
A combined model based on the extended finite element method, the mixed-mode cohesive zone model and the equivalent single layer theory is proposed in \cite{shells_fem_combine}. This model allows to obtain both initiation and propagation of delaminations. The novelty in the approach is that the intact material is modeled with a flat-shell formulation based on a first-order shear deformation theory, and elements that satisfy a certain criterion for de-bonding are switched to a more complex and computationally expensive extended finite element method with enhanced degrees of freedom that represents the delaminated material.

The review in \cite{taha_2020} discusses the historical development of shell theory, from Kirchhoff-Love models to more contemporary approaches. One direction discussed is the formulation of shell models as hyperbolic systems of equations that allows to separate different types of waves mathematically. In \cite{semenov_1990} a method is proposed to reduce a Timoshenko-type shell model to a first-order hyperbolic system. The problem with this approach is that it considers an axisymmetric statement, and its extension to a full statement greatly increases the number of variables. 

In this paper, we propose a method for reducing the Kirchhoff--Love equations \cite{kirchhoff_1850,love_1888} to a first-order hyperbolic system, and an approach to solving this system using the grid-characteristic method \cite{vas_2016}. The advantage of our approach is its low computational complexity. The aim of this study is to compare the accuracy of our simplified reduction with a full solution obtained from three-dimensional elasticity theory.

\section{Mathematical model}
\subsection{Three-dimensional elasticity}

The full description of the adduction of the viscoelasticity model to a hyperbolic system of equations is described in \cite{vas_2016}.
The dynamic equations for a linear elasticity model:
\begin{eqnarray}
\label{initial_equations}
\rho\dot{v}_i &= \nabla_j\sigma_{ij} & \textrm{(motion equations)}\nonumber\\
\dot\sigma_{ij} &= q_{ijkl}\dot{\varepsilon}_{kl} & \textrm{(rheological equations).}
\end{eqnarray}

In these equations 
$\rho$ is the medium density, 
$v_i$ are the displacement velocity components,
$\sigma_{ij}$, $\varepsilon_{ij}$ are the components of the stress tensors and the strain velocities,
$\nabla_j$ is the covariance derivative with respect to the $j$th coordinate, 
tensor $q_{ijkl}$ determines the rheology of the medium.

For linear elasticity, the tensor $q_{ijkl}$ in (\ref{initial_equations}) has the following form:

\begin{equation}
\label{tensor_qijkl_elastic}
q_{ijkl} = \lambda\delta_{ij}\delta_{kl}+\mu(\delta_{ik}\delta_{jl}+\delta_{il}
\delta_{jk}),\nonumber\\
\end{equation}
Here $\lambda$ and $\mu$ are Lame parameters, and $\delta_{ij}$ is Kronecker symbol.

The key assumption that allows to reduce the system (\ref{initial_equations}) to a hyperbolic form is as follows.
For small strains, tensor $e_{ij}=\dot{\varepsilon}_{ij}$ 
can be expressed in the following form:
\begin{equation}
e_{ij}=\frac{1}{2}(\nabla_j v_i+\nabla_i v_j).
\end{equation}

This substitution allows to reduce the vector of variables to nine components: $\vec{u} = (v_1, v_2, v_3, \sigma_{11}, \sigma_{12}, \sigma_{13}, \sigma_{22}, \sigma_{23}, \sigma_{33})^T$ and rewrite the original system (\ref{initial_equations}) in the matrix form:

\begin{equation}
\label{matrix_equation}
\frac{\partial\vec{u}}{\partial{t}}+\mathbf{A}_x\frac{\partial\vec{u}}{\partial{x}}+
\mathbf{A}_y\frac{\partial\vec{u}}{\partial{y}}+
\mathbf{A}_z\frac{\partial\vec{u}}{\partial{z}}=\vec{f}.
\end{equation}
In this form $\mathbf{A}_x$, $\mathbf{A}_y$, $\mathbf{A}_z$ are the matrices of ninth order that are constant and depend only on material properties.

\subsection{Kirchhoff-Love thin plates theory}

The theory, developed by A. E. H. Love \cite{love_1888}, is a two-dimensional elasticity model for thin plates and relies on the following postulates proposed by G. Kirchhoff \cite{kirchhoff_1850}:

\begin{itemize}
    \item lines normal to the middle surface remain straight and perpendicular to the deformed middle surface;
    \item transverse shear strains and the normal stress in the thickness direction are zero; consequently, the plate thickness does not change during deformation.
\end{itemize}

This model is not sufficiently rigorous for large deformations and thick plates, but is well-suited for describing small deflections and the propagation of elastic waves.

The system of equations describing the relationship between stress resultants and strains, and between moments and curvatures, is given by \cite{timoshenko_1959}:

\begin{eqnarray}
N_x &=& \frac{Eh}{1 - \nu^2} (\varepsilon_1 + \nu \varepsilon_2),
\label{eq:reology1}\\
N_y &=& \frac{Eh}{1 - \nu^2} (\varepsilon_2 + \nu \varepsilon_1),
\label{eq:reology2}\\
N_{xy} &=& N_{yx} = \frac{\gamma h E}{2(1 + \nu)}, 
\label{eq:reology3}\\
M_x &=& -D (\chi_x + \nu \chi_y),
\label{eq:reology4}\\
M_y &=& -D (\chi_y + \nu \chi_x),
\label{eq:reology5}\\
M_{xy} &=& -M_{yx} = D(1 - \nu) \chi_{xy},
\label{eq:reology6}
\end{eqnarray}

where:
\begin{itemize}
    \item \(N_x, N_y\) -- resultant forces per unit length on normal sections in the x and y directions, respectively (normal stress resultants).
    \item \(N_{xy}, N_{yx}\) -- shearing (tangential) forces in the midplane per unit length (shear stress resultants).
    \item \(M_x, M_y\) -- bending moments per unit length acting on faces perpendicular to the x and y axes, respectively.
    \item \(M_{xy}, M_{yx}\) -- twisting moments per unit length.
    \item \(E\) -- elastic modulus (Young’s modulus) of the material.
    \item \(\nu\) -- Poisson’s ratio of the material.
    \item \(h\) -- plate thickness.
    \item \(\varepsilon_1, \varepsilon_2\) -- normal strains of the midplane in the x and y directions.
    \item \(\gamma\) -- angular deformation (shear strain) in the midplane in the xy coordinate system. In the equation for \(N_{xy}\), the quantity \(\frac{E}{2(1 + \nu)}\)  represents the shear modulus \(G\).
    \item \(D\) -- flexural rigidity (bending stiffness) of the plate, defined as \(D = \frac{Eh^3}{12(1-\nu^2)}\).
    \item \(\chi_x, \chi_y\) -- changes in curvatures of the middle surface in the x and y directions.
    \item \(\chi_{xy}\) -- change in twisting curvature of the middle surface.
\end{itemize}

Let us define the stress tensor $\mathbf{T}$ and the moment tensor $\mathbf{M}$ as follows:

\begin{center}
$$
    \mathbf{T} = \left( \begin{array}{cc}
              \sigma_x & \sigma_{xy} \\
              \sigma_{xy} & \sigma_y
         \end{array} \right)
= \frac{1}{h} \left( \begin{array}{cc}
                         N_x & N_{xy} \\
                         N_{xy} & N_y
                      \end{array} \right).
$$
\end{center}

\begin{center}
$$
    \mathbf{M} = \left( \begin{array}{cc}
              M_x & M_{xy} \\
              M_{yx} & M_y
         \end{array} \right).
$$
\end{center}

Therefore, the rheological relations (\ref{eq:reology1} - \ref{eq:reology6}) can be written in the form:

\begin{equation}
\label{eq:T}
    \mathbf{T} = E \left( \begin{array}{cc} 
                 \frac{\varepsilon_1 + \nu \varepsilon_2}{1 - \nu^2} & \frac{\gamma}{2(1 + \nu)} \\ \\
                 \frac{\gamma}{2(1 + \nu)} & \frac{\varepsilon_2 + \nu \varepsilon_1}{1 - \nu^2}
                 \end{array} \right),
\end{equation}
\\
\begin{equation}
\label{eq:M}
    \mathbf{M} = D \left( \begin{array}{cc}
                 -\chi_x - \nu \chi_y & (1 - \nu) \chi_{xy} \\ \\
                 (1 - \nu) \chi_{xy} & -\chi_y - \nu \chi_x
                 \end{array} \right).
\end{equation}

The dynamic relations can be written as:

\begin{eqnarray}
\rho \dot{\vec{v}} &=& \nabla \cdot \mathbf{T} \label{eq:vel}\\
I \dot{\vec{\omega}} &=& \nabla \cdot \mathbf{M} \label{eq:ang_vel},
\end{eqnarray}

where:
\begin{itemize}
    \item $\rho$ - material density;
    \item $I = \frac{\rho h^3}{12}$ - moment of inertia of the unit area of the plate about its midplane;
\end{itemize}

In order to bring all equations to a unified form and write them in matrix form, we use the relation for the small strain tensor. Its time derivative can be written as \cite{chelnokov_2005}:

\begin{equation}
\label{eq:strain_tensor}
    \dot{\mathbf{\varepsilon}} = \left( \begin{array}{cc}
    \frac{\partial v_x}{\partial x} & \frac{1}{2} \left( \frac{\partial v_x}{\partial y} + \frac{\partial v_y}{\partial x} \right) \\ \\
    \frac{1}{2} \left( \frac{\partial v_x}{\partial y} + \frac{\partial v_y}{\partial x} \right) & \frac{\partial v_y}{\partial y}
    \end{array} \right)
\end{equation}

where $v_x, v_y$ are the components of the velocity vector.

The quantities $\chi$ can be similarly written in the form of a curvature tensor, whose time derivative is:

\begin{equation}
\label{eq:curv_tensor}
    \dot{\mathbf{\chi}} = \left( \begin{array}{cc}
    \frac{\partial \omega_x}{\partial x} & \frac{1}{2} \left( \frac{\partial \omega_x}{\partial y} + \frac{\partial \omega_y}{\partial x} \right) \\ \\
    \frac{1}{2} \left( \frac{\partial \omega_x}{\partial y} + \frac{\partial \omega_y}{\partial x} \right) & \frac{\partial \omega_y}{\partial y}
        \end{array} \right)
\end{equation}

where $w_x, w_y$ are the components of the angular velocity vector.

Using relations (\ref{eq:strain_tensor}) and (\ref{eq:curv_tensor}), the expressions for the derivatives of the stress and moment tensors (\ref{eq:T}) and (\ref{eq:M}) can be written as follows:

\begin{equation}
    \dot{\mathbf{T}} = E \left( \begin{array}{cc} 
                 \frac{\frac{\partial v_x}{\partial x} + \nu \frac{\partial v_y}{\partial y}}{1 - \nu^2} 
                 & 
                 \frac{\frac{\partial v_y}{\partial x} + \frac{\partial v_x}{\partial y}}{4(1 + \nu)} \\ \\
                 \frac{\frac{\partial v_y}{\partial x} + \frac{\partial v_x}{\partial y}}{4(1 + \nu)} 
                 & 
                 \frac{\nu  \frac{\partial v_x}{\partial x} + \frac{\partial v_y}{\partial y}}{1 - \nu^2}
                 \end{array} \right) 
                 \label{eq:dT_dt} \\ \\
\end{equation}

\begin{equation}
\label{eq:dM_dt}
\dot{\mathbf{M}} = D \left( \begin{array}{cc} 
                 - \frac{\partial w_x}{\partial x} - \nu \frac{\partial w_y}{\partial y} 
                 & 
                 \frac{1-\nu}{2} \left( \frac{\partial w_y}{\partial x} + \frac{\partial w_x}{\partial y} \right) \\ \\
                 \frac{\nu-1}{2} \left( \frac{\partial w_y}{\partial x} + \frac{\partial w_x}{\partial y} \right) 
                 & 
                 - \frac{\partial w_y}{\partial y} - \nu \frac{\partial w_x}{\partial x} 
                 \end{array} \right)
\end{equation}

Note that (\ref{eq:vel}), (\ref{eq:ang_vel}), (\ref{eq:dT_dt}), (\ref{eq:dM_dt}) all have the same shape, which allows us to write them in a matrix form:

\begin{equation}
\label{eq:hyperbolic}
    \frac{\partial \vec{U}}{\partial t} + A_x \frac{\partial \vec{U}}{\partial x} + A_y \frac{\partial \vec{U}}{\partial y} = \vec{0}
\end{equation}

where:

\begin{itemize}
    \item $\vec{U}$ = $\left( v_x, v_y, w_x, w_y, \sigma_x, \sigma_y, \sigma_{xy}, M_x, M_y, M_{xy}\right )$;
    \item $A_x$ and $A_y$ are written below.
\end{itemize}

\begin{equation}
\label{eq:AX}
A_x = -
\left(
\begin{array}{cccccccccc}
\cdot & \cdot & \cdot & \cdot & \frac{1}{\rho} & \cdot & \cdot & \cdot & \cdot & \cdot \\
\cdot & \cdot & \cdot & \cdot & \cdot & \cdot & \frac{1}{\rho} & \cdot & \cdot & \cdot \\
\cdot & \cdot & \cdot & \cdot & \cdot & \cdot & \cdot & \frac{1}{I} & \cdot & \cdot \\
\cdot & \cdot & \cdot & \cdot & \cdot & \cdot & \cdot & \cdot & \cdot & \frac{1}{I} \\
\frac{E}{1-\nu^2} & \cdot & \cdot & \cdot & \cdot & \cdot & \cdot & \cdot & \cdot & \cdot \\
\frac{E\nu}{1-\nu^2} & \cdot & \cdot & \cdot & \cdot & \cdot & \cdot & \cdot & \cdot & \cdot \\
\cdot & \frac{E}{4(1+\nu)} & \cdot & \cdot & \cdot & \cdot & \cdot & \cdot & \cdot & \cdot \\
\cdot & \cdot & D & \cdot & \cdot & \cdot & \cdot & \cdot & \cdot & \cdot \\
\cdot & \cdot & D\nu & \cdot & \cdot & \cdot & \cdot & \cdot & \cdot & \cdot \\
\cdot & \cdot & \cdot & \frac{D}{2}(1-\nu) & \cdot & \cdot & \cdot & \cdot & \cdot & \cdot \\
\end{array}
\right)\end
{equation}

\begin{equation}
\label{eq:AY}
A_y = -
\left(
\begin{array}{cccccccccc}
\cdot & \cdot & \cdot & \cdot & \cdot & \cdot & \frac{1}{\rho} & \cdot & \cdot & \cdot \\
\cdot & \cdot & \cdot & \cdot & \cdot & \frac{1}{\rho} & \cdot & \cdot & \cdot & \cdot \\
\cdot & \cdot & \cdot & \cdot & \cdot & \cdot & \cdot & \cdot & \cdot & \frac{1}{I} \\
\cdot & \cdot & \cdot & \cdot & \cdot & \cdot & \cdot & \cdot & \frac{1}{I} & \cdot \\
\cdot & \frac{E\nu}{1-\nu^2} & \cdot & \cdot & \cdot & \cdot & \cdot & \cdot & \cdot & \cdot \\
\cdot & \frac{E}{1-\nu^2} & \cdot & \cdot & \cdot & \cdot & \cdot & \cdot & \cdot & \cdot \\
\frac{E}{4(1+\nu)} & \cdot & \cdot & \cdot & \cdot & \cdot & \cdot & \cdot & \cdot & \cdot \\
\cdot & \cdot & \cdot & D\nu & \cdot & \cdot & \cdot & \cdot & \cdot & \cdot \\
\cdot & \cdot & \cdot & D & \cdot & \cdot & \cdot & \cdot & \cdot & \cdot \\
\cdot & \cdot & \frac{D}{2}(1-\nu) & \cdot & \cdot & \cdot & \cdot & \cdot & \cdot & \cdot \\
\end{array}
\right)\end
{equation}

These matrices have basis from eigen vectors and all their eigen values are real, so this is a hyperbolycal system and can be solved via GCM. 

\section{Numerical method}

The grid-characteristic numerical method (GCM) is described in \cite{chelnokov_2005}, and its formulation for irregular tetrahedral grids is described in \cite{vas_2016}.
We use irregular grids in this research to eliminate numerical anisotropy.

The implementation of GCM for shells uses the same procedure as in \cite{vas_2016}, but with material martices (\ref{eq:AX}) and (\ref{eq:AY}).

\subsection{Moments calculation on the interpolation grid}

For comparison between 3d and shells moment waves images required to calculate the first one on non-regular grid approximated plate.

From \cite{timoshenko_1959} we can take the following formulas:
$$M_{xx} = \int_{-h/2}^{h/2}\sigma_{xx}z(1-\frac{z}{r_y})dz, \ M_{yy} = \int_{-h/2}^{h/2}\sigma_{yy}z(1-\frac{z}{r_x})dz$$
$$M_{xy} = -\int_{-h/2}^{h/2}\tau_{xy}z(1-\frac{z}{r_y})dz, \  M_{yx} = \int_{-h/2}^{h/2}\tau_{yx}z(1-\frac{z}{r_x})dz$$

In approximate that $h \ll r$ :
$$M_{xx} = \int_{-h/2}^{h/2}\sigma_{xx}zdz, \ M_{yy} = \int_{-h/2}^{h/2}\sigma_{yy}zdz$$
$$M_{xy} = -\int_{-h/2}^{h/2}\tau_{xy}zdz, \  M_{yx} = \int_{-h/2}^{h/2}\tau_{yx}zdz$$

Let's introduce the notation:
$$\sigma(h/2) = s_{up} ,\ \sigma(-h/2) = s_{down}$$

Let approximate a two-dimensional stress diagram along through lines drawn along the normals at each point of median plane using a linear function. Then it will be possible to directly integrate the stresses along these lines and obtain a moment at each point of this plane:

   $$\sigma(z) = \frac{s_{up} - s_{down}}{h} z + \frac{s_{up} + s_{down}}{2}$$
$$M = \int_{-h/2}^{h/2}\sigma zdz = \int\limits_{-h/2}^{h/2} f(z) \cdot z \, dz = \frac{(s_{\mathrm{up}} - s_{\mathrm{down}})\, h^2}{12}$$

It is worth saying that such a calculation of moments becomes correct only in the case of sufficiently thin plates, what will be shown next.

\section{Numerical experiments}

The hyperbolic set of equations (\ref{eq:hyperbolic}) separates into two independent parts.

The first one describes the dependency between velocities and stresses. 
It can be associated with a two dimensional model and "in-plane" elastic waves that propagate inside the plate.

The second part of the set ties angular velocities and moments. 
This part describes "momentum", "out-of-plane" waves that can be associated with Lamb waves that propagate with the same speed.

It is worth noting that in a flat plate these parts are independent, so in-plane and momentum waves propagate indepentently.

To verify this separation, two different statements were considered and shell calculations were compared to plates of different thickness calculated with a full three dimensional method.

\subsection{In-plane disturbance}
In this section the models were compared in the following setup:

\begin{enumerate}
    \item a thin plate of size $10\times 10\ \mathrm{m}$ (thickness varied for the three-dimensional computation);
    \item material parameters: density $\rho = 7800\ \mathrm{kg\,m^{-3}}$, Young's modulus $E = 210\ \mathrm{GPa}$, Poisson's ratio $\nu = 0.30$. The material is isotropic;
    \item initial condition: a velocity $V_x = 100\ \mathrm{m\,s^{-1}}$ is prescribed at the central point of the plate, directed along the $x$-axis;
    \item boundary conditions: zero-gradient (Neumann) boundary conditions.
\end{enumerate}

\subsubsection{Three-dimensional model}

When discussing complete three-dimensional wave patterns in thin plates, our main task is to identify those features that cannot be taken into account by a two-dimensional shell-model and lead to discrepancies.

Consider the effect of the "tail" and the stretching of the wave front. To do this, let us look at the wave pattern in the cross sections of the plate depending on its thickness. Consider three different thicknesses of 1.0, 0.7 and 0.4 meters on three different moments in time that demonstrate the effect.

\begin{figure}[!h]
    \centering
    \includegraphics[width=\linewidth]{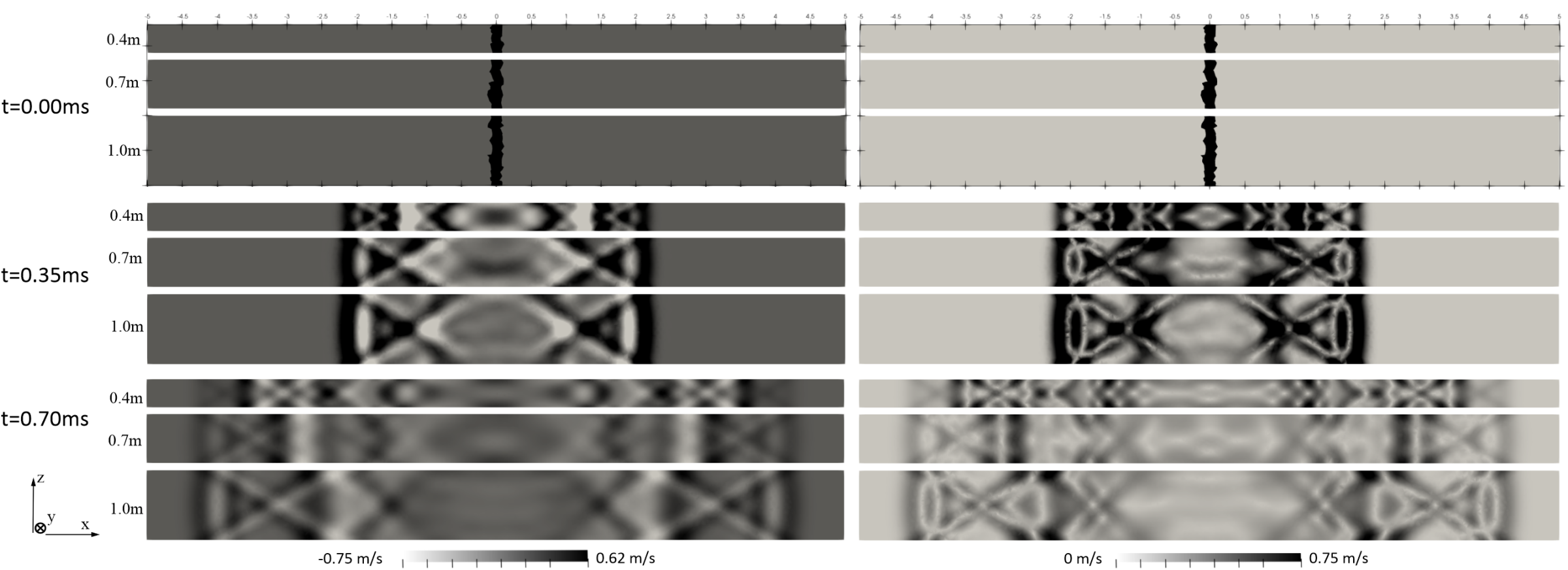}
    \caption{Wave pattern on a XZ slice through the impact point: on the left - the x-axis component of velocity, on the right - magnitude of velocity}
    \label{fig:inplane_ynorm_all}
\end{figure}

This effect can be explained by the fact that in a homogeneous isotropic medium in the absence of boundaries, the front of shear and longitudinal elastic waves tends to propagate spherically from a point source. Multiple reflections from the upper and lower boundaries of the plate lead to displacement of the maximum amplitude of the wave relative to its leading edge, forming a "tail", which can be interpreted as a stretched and displaced wave front that propagates slower than the three dimensional wave. \\

This effect turned out to be less noticeable for the transverse wave. We can see this by examining its propagation in the second cross-section of the plate.

\newpage

\begin{figure}[!h]
    \centering
    \includegraphics[width=0.8\linewidth]{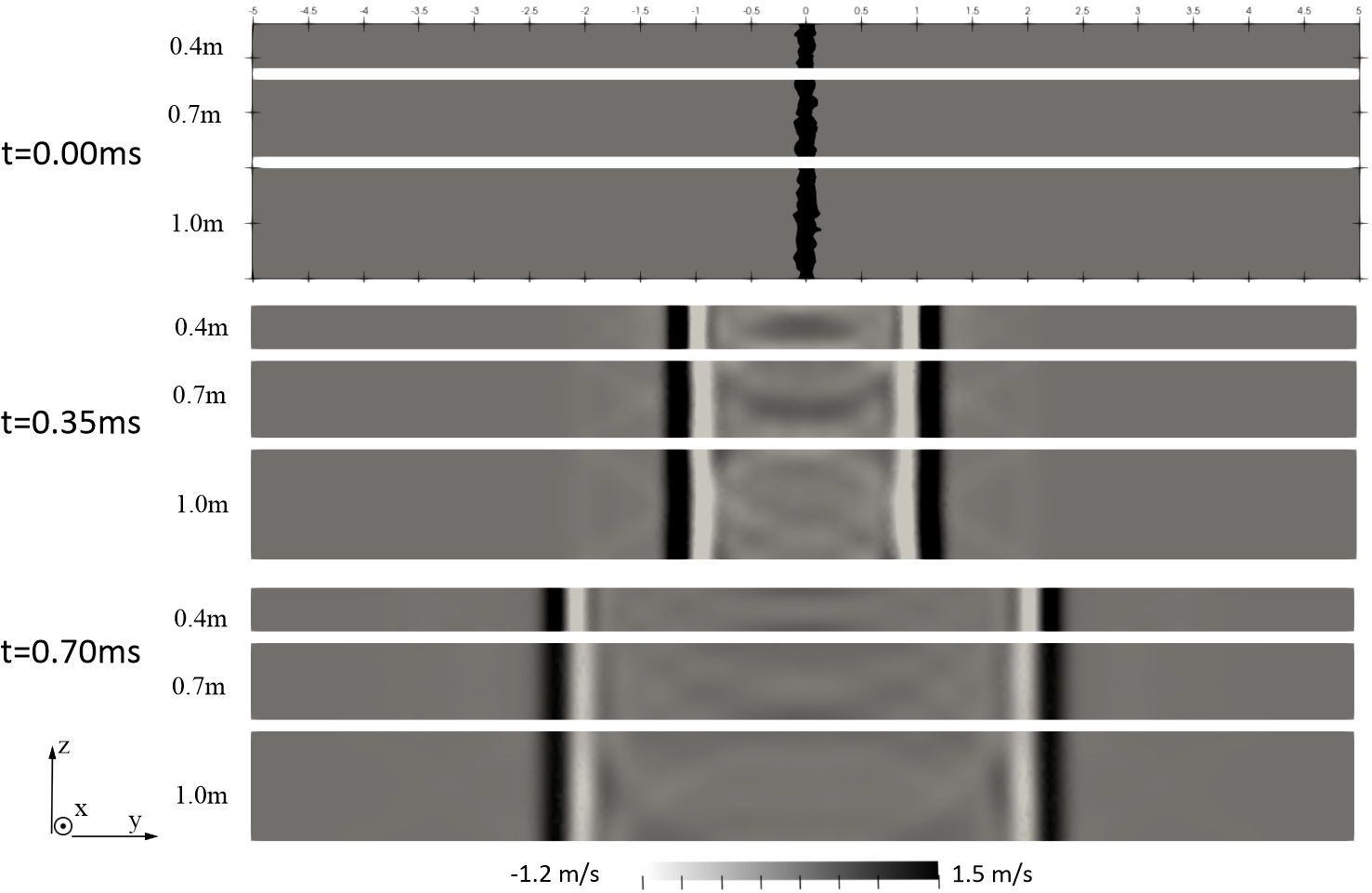}
    \caption{Wave pattern on a YZ slice through the impact point. X-axis component of velocity.}
    \label{fig:inplane_xnorm_vx_all}
\end{figure}

\begin{figure}[!h]
    \centering
    \includegraphics[width=\linewidth]{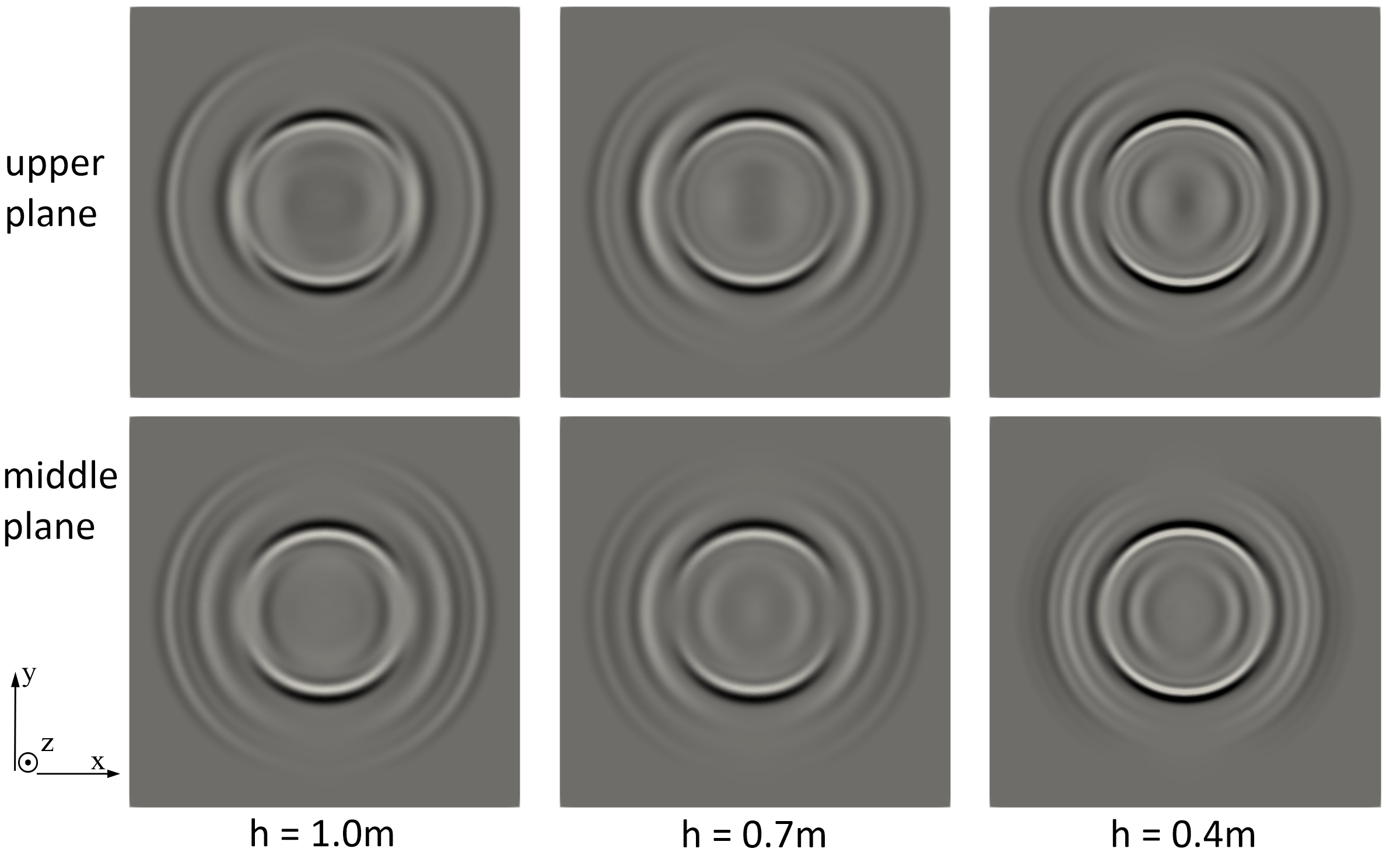}
    \caption{Distribution of X-axis component of velocity in different sections at 0.70 ms}
    \label{fig:inplane_upv}
\end{figure}

Figure \ref{fig:inplane_upv} shows the general view of the wave pattern at the maximum time in the pictures above. We will display the view of the upper and middle sections of the plate, but the results of the shell model are compared only with the second one due to the basic formulation of the model.

\subsubsection{Kirchhoff–Love model}

\begin{figure}[!h]
    \centering
    \begin{minipage}{0.3\linewidth}
        \centering
        \includegraphics[width=\linewidth]{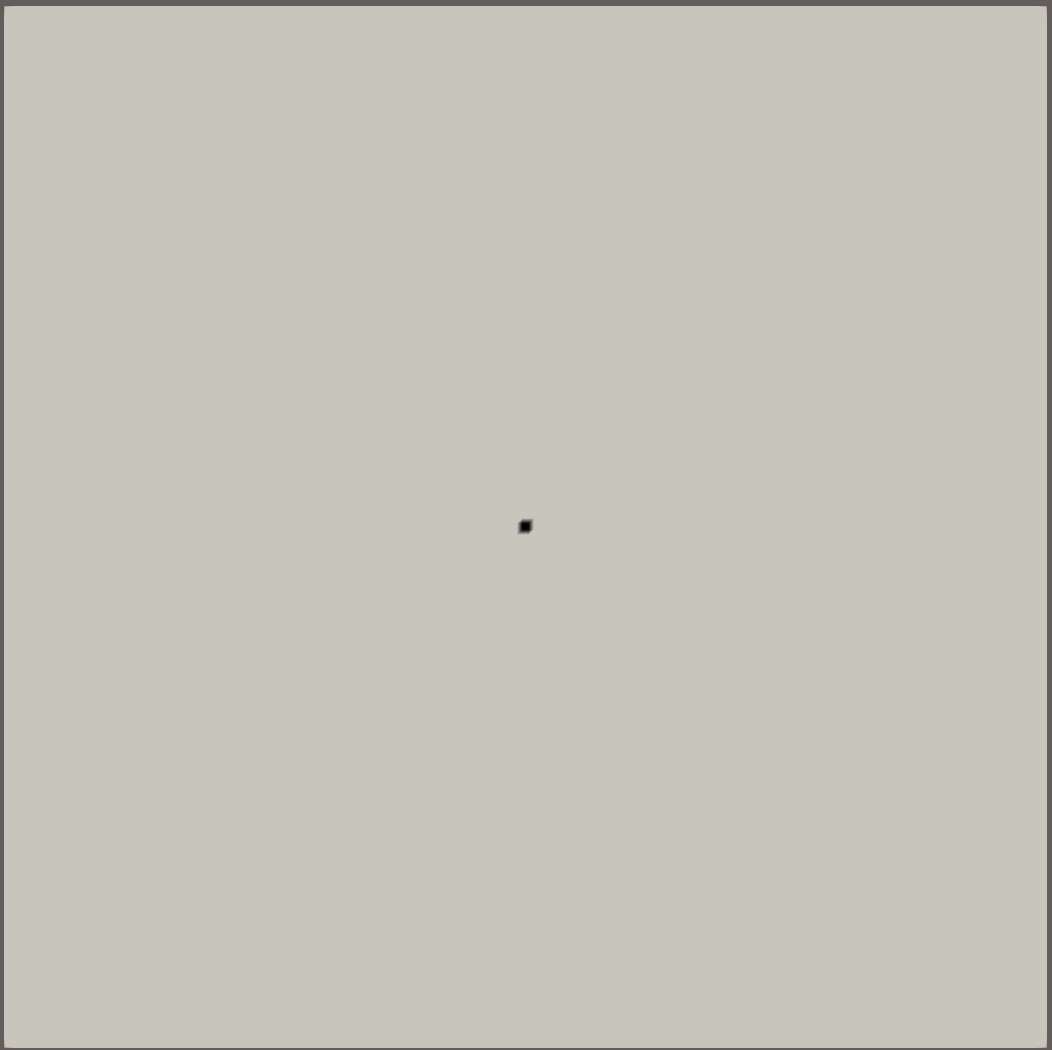}
        \caption{Velocity magnitude $t = 0\ \mathrm{s}$}
        \label{fig:kl_Vx_0}
    \end{minipage}
    \hfill
    \begin{minipage}{0.3\linewidth}
        \centering
        \includegraphics[width=\linewidth]{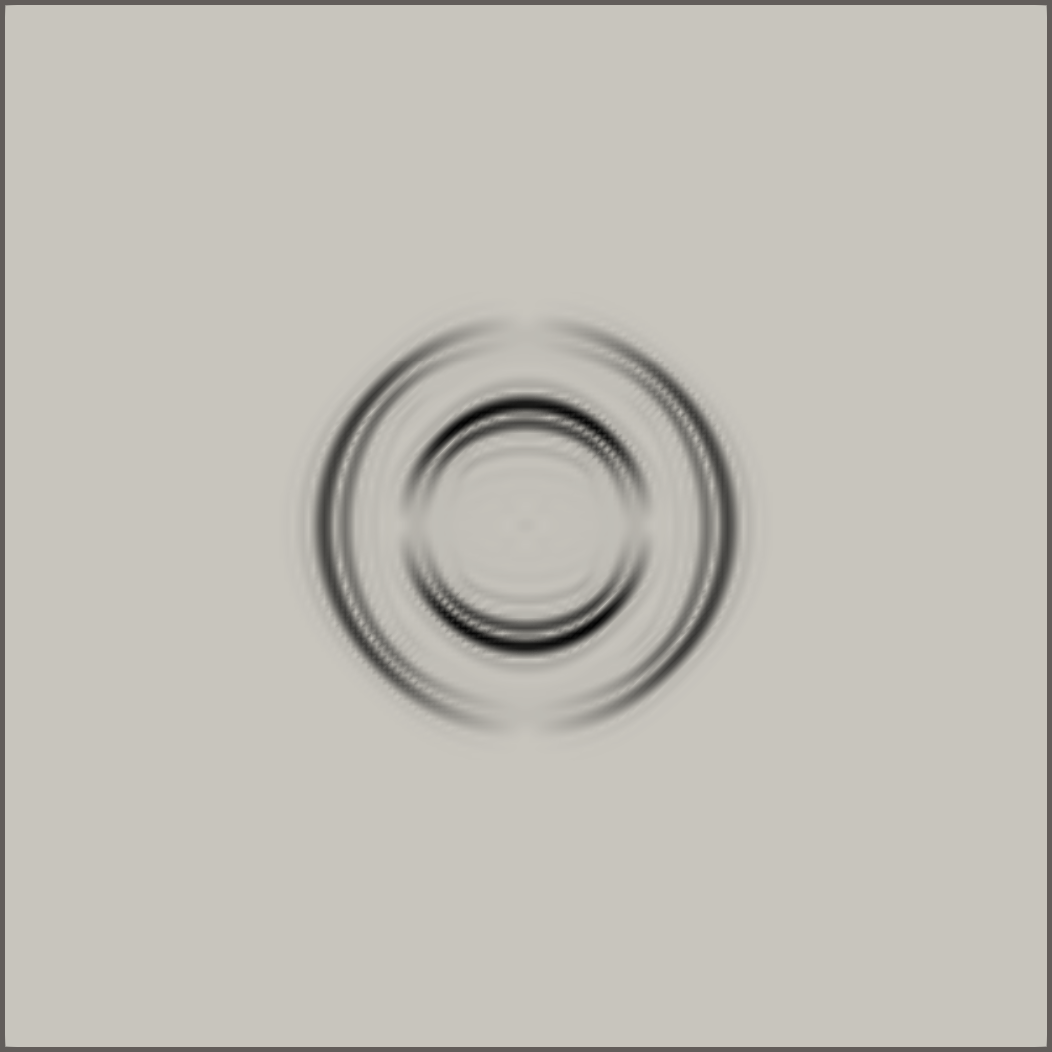}
        \caption{Velocity magnitude  $t = 0.35 \mathrm{ms}$}
        \label{fig:kl_Vx_175}
    \end{minipage}
    \hfill
    \begin{minipage}{0.3\linewidth}
        \centering
        \includegraphics[width=\linewidth]{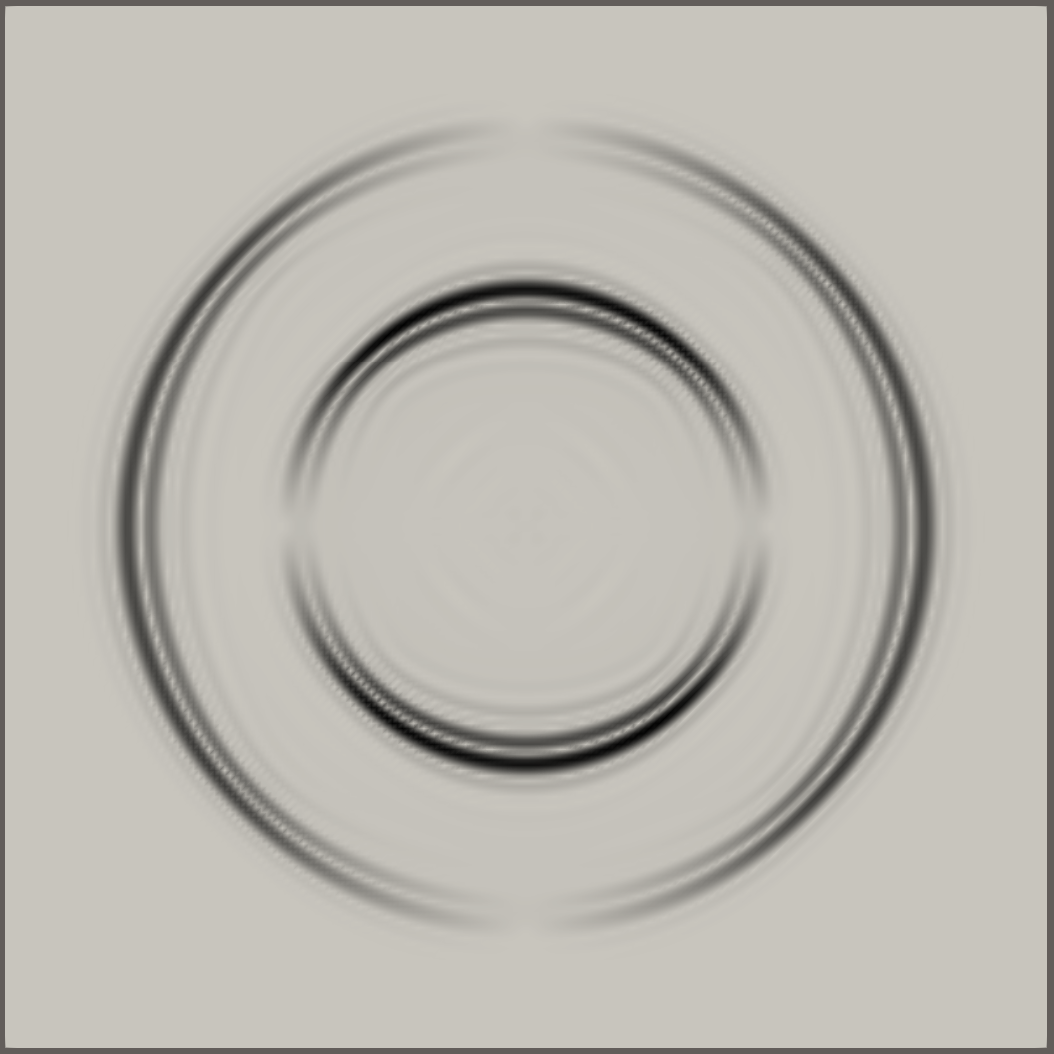}
        \caption{Velocity magnitude  $t = 0.7 \mathrm{ms}$}
        \label{fig:kl_Vx_350}
    \end{minipage}
\end{figure}

\begin{figure}[!h]
    \centering
    \begin{minipage}{0.3\linewidth}
        \centering
        \includegraphics[width=\linewidth]{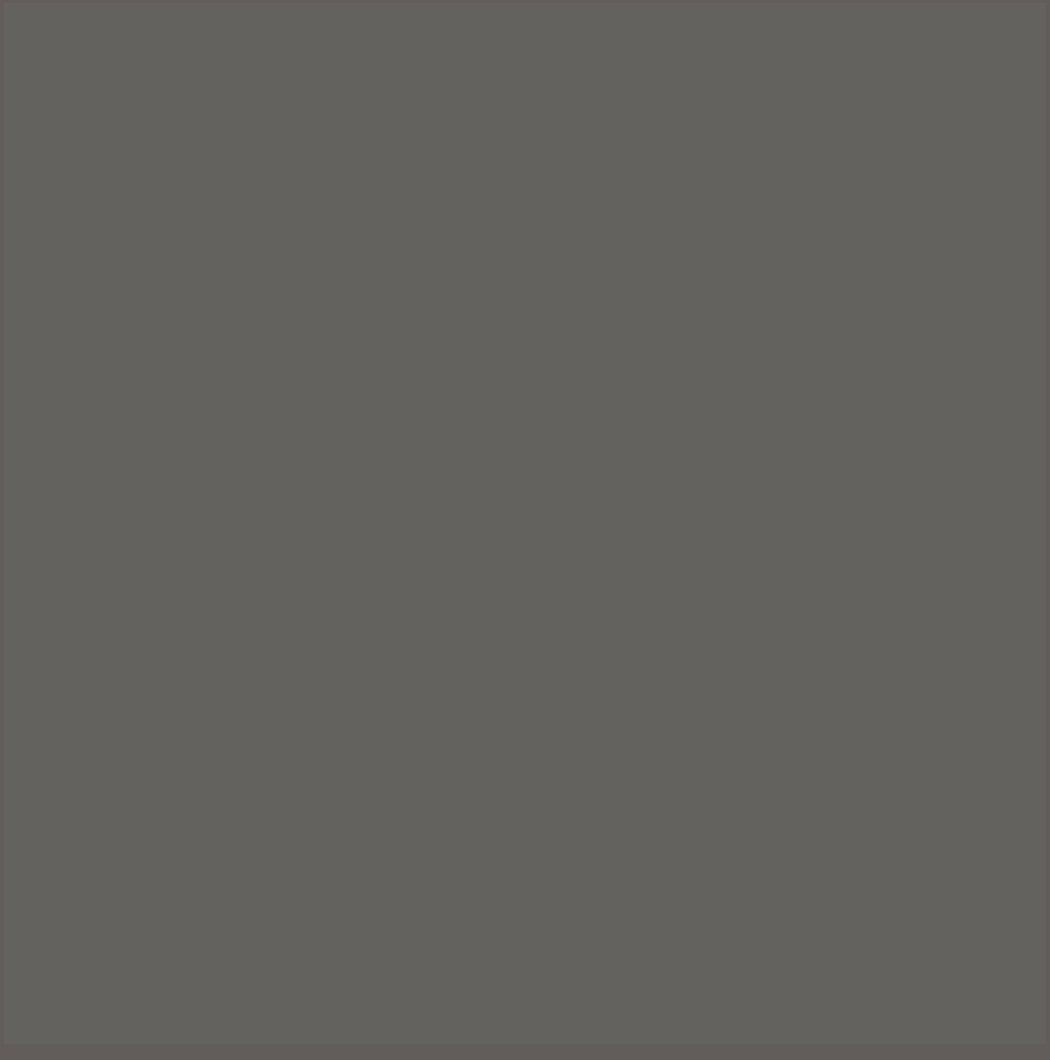}
        \caption{$M_{xx}$ $t = 0\ \mathrm{s}$}
        \label{fig:kl_Mxx_0}
    \end{minipage}
    \hfill
    \begin{minipage}{0.3\linewidth}
        \centering
        \includegraphics[width=\linewidth]{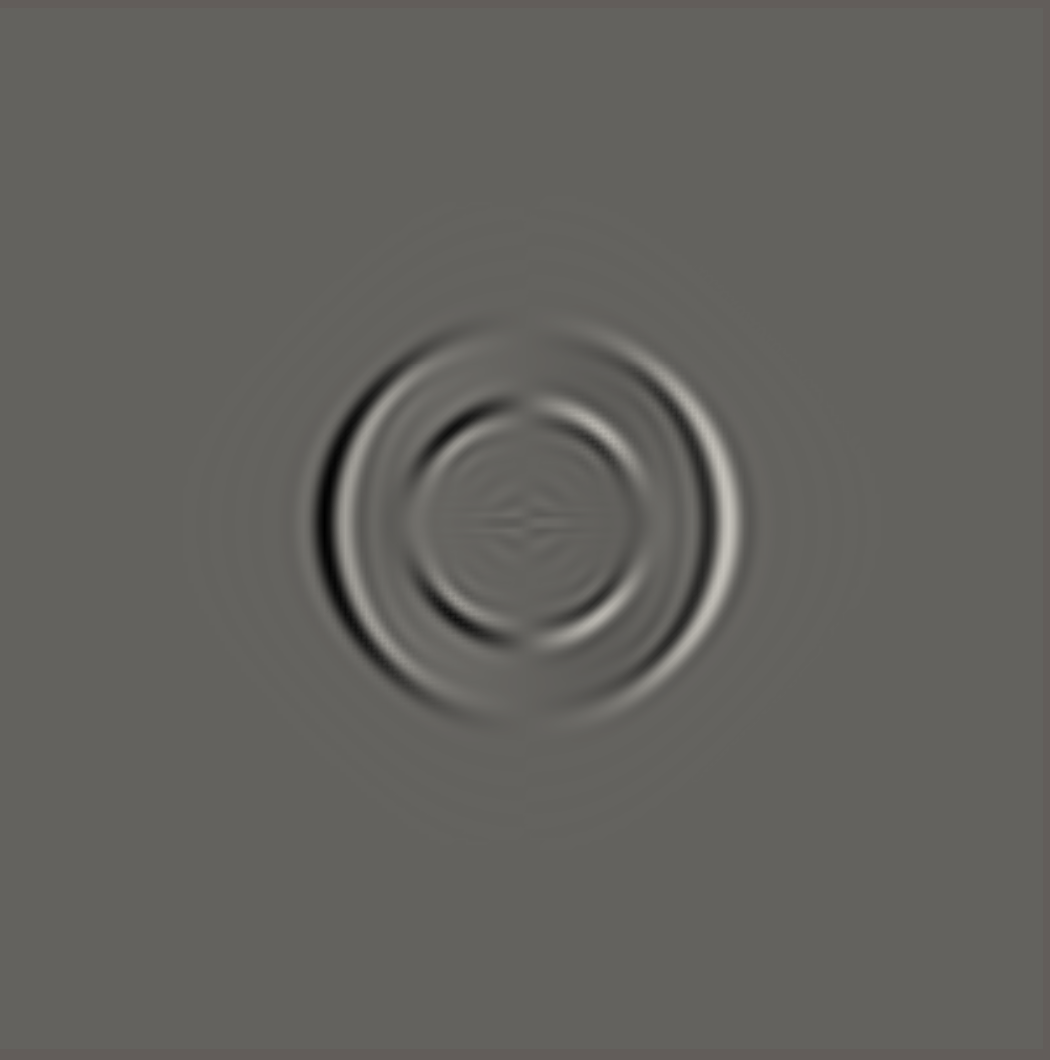}
        \caption{$M_{xx}$ $t = 0.35 \mathrm{ms}$}
        \label{fig:kl_Mxx_175}
    \end{minipage}
    \hfill
    \begin{minipage}{0.3\linewidth}
        \centering
        \includegraphics[width=\linewidth]{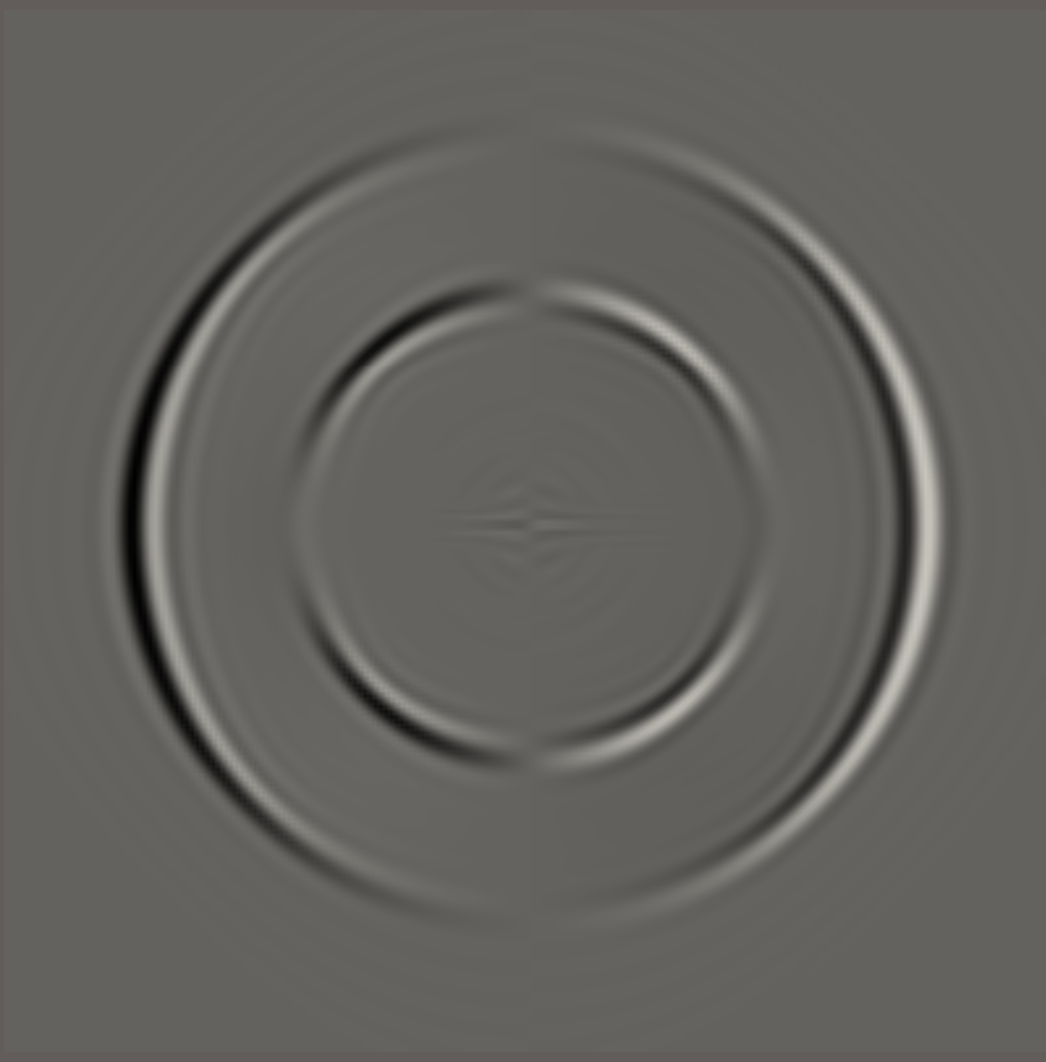}
        \caption{$M_{xx}$ $t = 0.7 \mathrm{ms}$}
        \label{fig:kl_Mxx_350}
    \end{minipage}
\end{figure}

To solve the problem within the framework of the Kirchhoff–Love model, the grid-characteristic method with fifth-order Newton polynomial interpolation was employed. The computed velocity amplitude at several time steps in a thin plate of size $10\times 10\ \mathrm{m}$ is shown in Figures~\ref{fig:kl_Vx_0}--\ref{fig:kl_Vx_350}. The computed bending moment $M_{xx}$ at the same time steps is shown in Figures~\ref{fig:kl_Mxx_0}--\ref{fig:kl_Mxx_350}.

We can clearly observe longitudinal waves radiating in the plate plane, as well as transverse waves propagating across the plate.

The wave propagation velocities correspond to the eigenvalues of the system matrix and are equal to $5439.3\ \mathrm{m\,s^{-1}}$ for longitudinal waves and $3217.9\ \mathrm{m\,s^{-1}}$ for transverse waves, respectively.

\subsubsection{Comparison of models}

\begin{figure}[!h]
    \centering
    \begin{minipage}{0.45\linewidth}
        \centering
        \includegraphics[width=\linewidth]{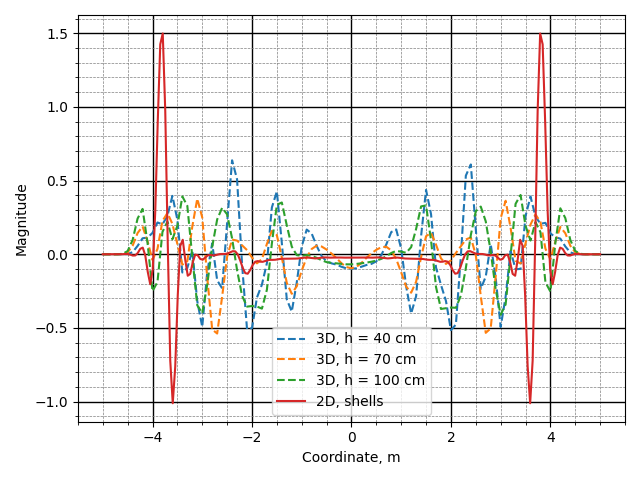}
        \caption{$V_x$ on X axis at $t = 0.7 \mathrm{ms}$}
        \label{fig:line_Vx_X}
    \end{minipage}
    \hfill
    \begin{minipage}{0.45\linewidth}
        \centering
        \includegraphics[width=\linewidth]{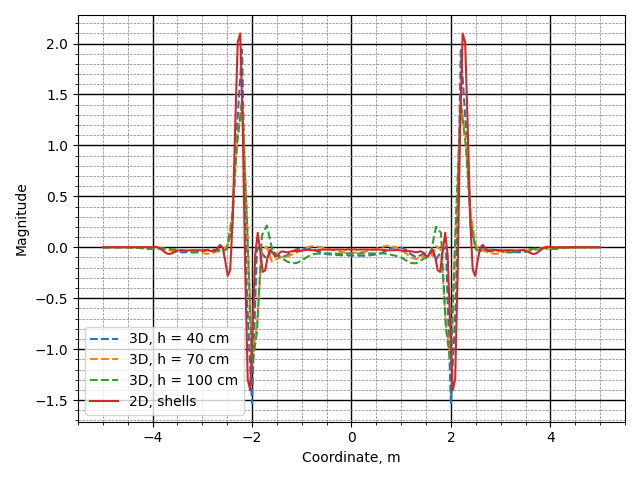}
        \caption{$V_x$ on Y axis at $t = 0.7 \mathrm{ms}$}
        \label{fig:line_Vx_Y}
    \end{minipage}
\end{figure}

We compare cross-sectional profiles along the $X$- and $Y$-axes for the three-dimensional and two-dimensional models. Figures~\ref{fig:line_Vx_X} and \ref{fig:line_Vx_Y} show cross-sectional profiles of the velocity $V_x$ along the $X$- and $Y$-axes, respectively, at time $t = 0.7 \mathrm{ms}$; the profiles were averaged over a width of $1\ \mathrm{m}$ to reduce oscillations.

Longitudinal waves on Figure \ref{fig:line_Vx_X} propagate along the $X$-axis. The two-dimensional model exhibits a single high-amplitude peak, whereas the three-dimensional model displays a large number of oscillations. Nevertheless, the leading edge of the wave packet in the three-dimensional model is close to the high-amplitude oscillation of the two-dimensional model, indicating agreement in the wave propagation speed.

Shear waves on Figure \ref{fig:line_Vx_Y} propagate along the $Y$-axis, and a pronounced peak of maximum amplitude, which coincides across all models, is clearly visible.

\begin{figure}[!h]
    \centering
    \begin{minipage}{0.45\linewidth}
        \centering
        \includegraphics[width=\linewidth]{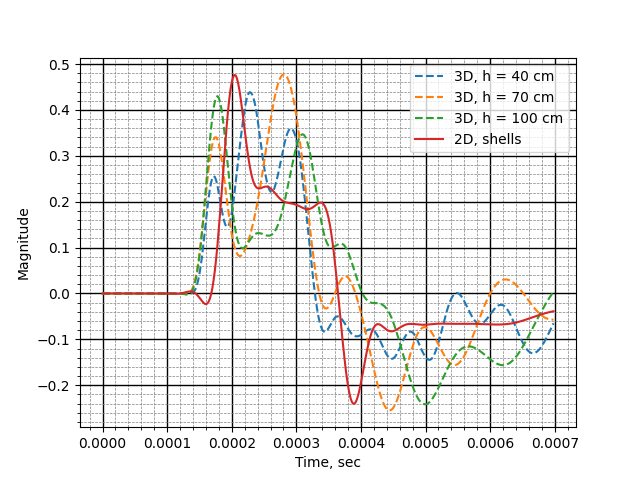}
        \caption{Sensor placed along X axis}
        \label{fig:sens_Vy_X}
    \end{minipage}
    \hfill
    \begin{minipage}{0.45\linewidth}
        \centering
        \includegraphics[width=\linewidth]{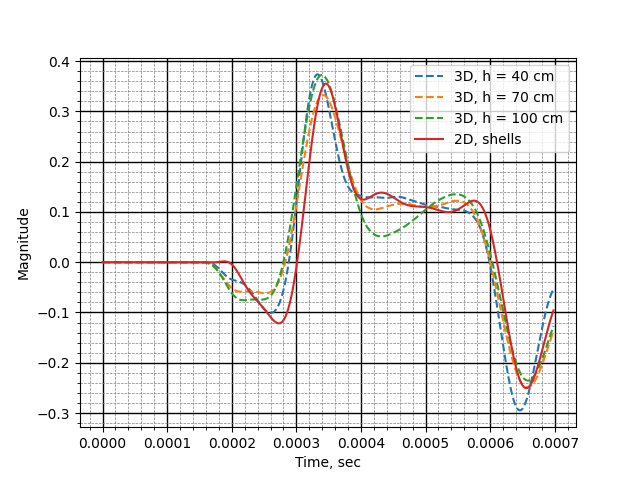}
        \caption{Sensor placed along Y axis}
        \label{fig:sens_Vy_Y}
    \end{minipage}
\end{figure}

In engineering applications, we rarely can have access to a full three dimensional pattern of stresses and velocities, so the comparison with an experiment is usually performed by writing a signal in a certain area of the material surface. In the present study two $1\times 1\ \mathrm{m}$ velocity sensors were modelled at two locations — translated by $(0,\,1.5)\ \mathrm{m}$ and $(1.5,\,0)\ \mathrm{m}$ from the plate centre where the initial condition was applied. Figures~\ref{fig:sens_Vy_X} and \ref{fig:sens_Vy_Y} show the time dependency of the spatially averaged amplitude $V_x$ for these sensor positions, respectively.

For the sensor located along the $Y$-axis the agreement between the two-dimensional model and the shell model is particularly notable. The sensor located along the $X$-axis recorded longitudinal waves which, as observed in the previous plots, exhibit strong oscillations.

As a quantitative measure of the agreement between the models we adopt the normalized root-mean-square error (NRMSE) \cite{zhou_nrmse}, which for two matrices can be expressed as follows:

\[
\mathrm{MSE} = \frac{1}{N_x N_y} \sum_{i=1}^{N_x} \sum_{j=1}^{N_y} 
\left( A_{ij} - B_{ij} \right)^2
\]

\[
\mathrm{RMSE} = \sqrt{\mathrm{MSE}}
\]

\[
\mathrm{NRMSE} = \frac{\mathrm{RMSE}}{\max\limits_{i,j} B_{ij} - \min\limits_{i,j} B_{ij}}
\]

Figure~\ref{fig:nrmse_Vx} shows the NRMSE metric as a function of time for different plate thicknesses in the three-dimensional simulations. It is evident that the error decreases as the thickness is reduced, which confirms that the Kirchhoff--Love model, as expected, describes shells more accurately the thinner they are.

\begin{figure}[!h]
    \centering
    \includegraphics[width=0.45\linewidth]{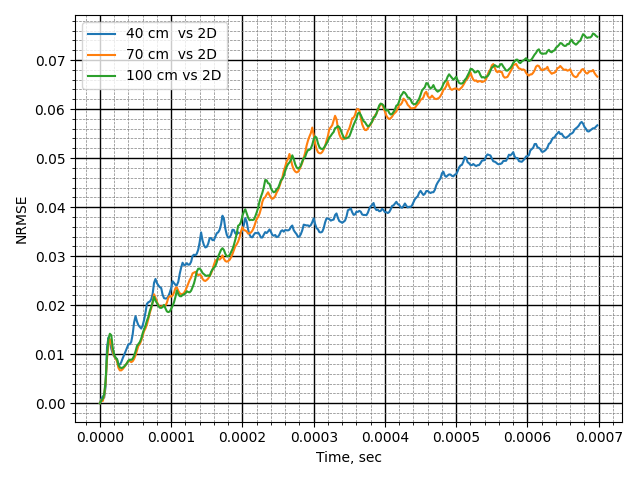}
    \caption{NRMSE metric for $V_x$}
    \label{fig:nrmse_Vx}
\end{figure}

\newpage
\subsection{Out-of-plane disturbance}

The statement is similar to the previous formulation, except that the initial velocity condition is now a linear gradient along the height from $100\ \mathrm{m\,s^{-1}}$ on the upper side and $-100\ \mathrm{m\,s^{-1}}$ on the bottom side.

\subsubsection{Three-dimensional model}

Let us now consider how the wave pattern evolves in the new statement. Here, the key for us is the reconciliation of bending and torques, which are integrals along the height of stresses, as shown in the paragraphs earlier. 

Consider the distribution of the normal component of the stress tensor sxx in the cross section of a plate with a normal along the y axis depending on thickness and time. There are no stresses on the first frame yet, so we will display the velocity on it to demonstrate the problem statement more clearly.

\begin{figure}[!h]
    \centering
    \includegraphics[width=\linewidth]{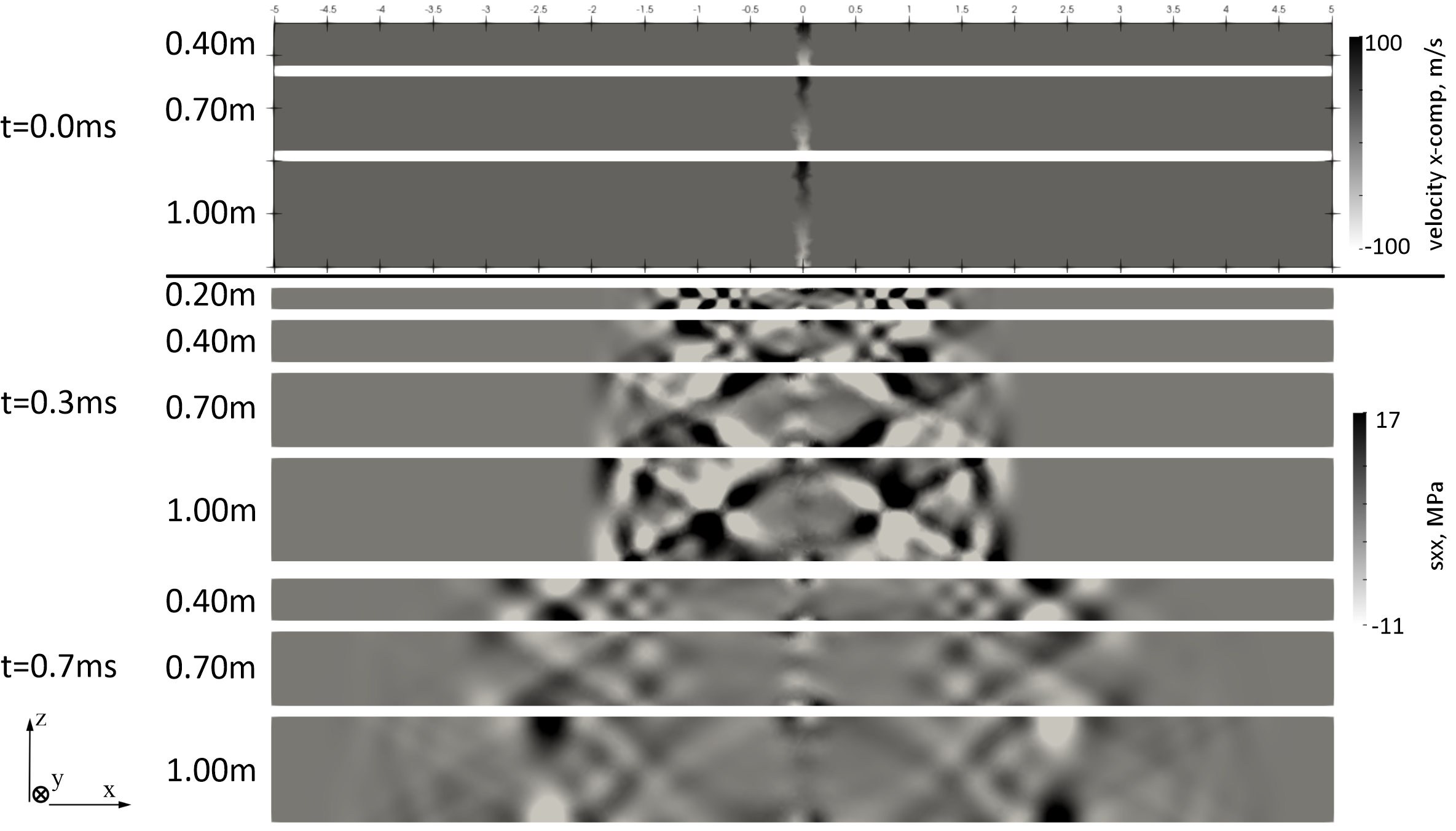}
    \caption{Distribution of sxx on a YZ slice through the impact point.}
    \label{fig:outplane_sxx_ynorm}
\end{figure}

Figure \ref{fig:outplane_sxx_ynorm} shows that when the plate is not thin enough, the stress distributions over height are strongly non-linear, which makes approximating them with a linear function for calculating moments incorrect. However, as the thickness decreases, we see a clear stabilization of the stress distribution in the section and the appearance of a kind of a mesh structure, which makes a linear approximation a correct model to to calculate the moments.

It is worth saying that the grid-like structure of the stress distribution is a moment wave, that is, with a decrease in the thickness of the plate, the moment waves stabilize. 

\newpage

\subsubsection{Kirchhoff-Love model}

From the standpoint of the Kirchhoff--Love model, angular velocities and bending moments are related to each other in the same manner as linear velocities and stresses. Therefore, the form of the angular velocity plots is analogous to that of the velocity plots. Let us now consider the time evolution of the bending moment $M_{xx}$.

\begin{figure}[!h]
    \centering
    \begin{minipage}{0.45\linewidth}
        \centering
        \includegraphics[width=\linewidth]{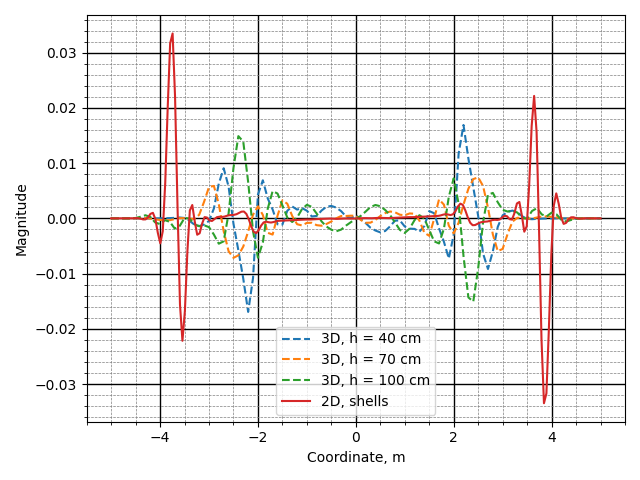}
        \caption{$M_{xx}$ on X axis at $t = 7.0\times10^{-4}\ \mathrm{s}$}
        \label{fig:line_Mxx_X}
    \end{minipage}
    \hfill
    \begin{minipage}{0.45\linewidth}
        \centering
        \includegraphics[width=\linewidth]{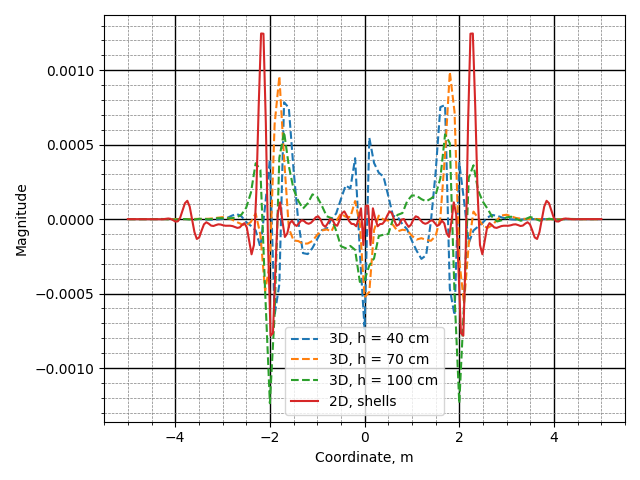}
        \caption{$M_{xx}$ on Y axis at $t = 7.0\times10^{-4}\ \mathrm{s}$}
        \label{fig:line_Mxx_Y}
    \end{minipage}
\end{figure}

Figures~\ref{fig:line_Mxx_X} and \ref{fig:line_Mxx_Y} show the bending moments along the coordinate axes at $t = 7.0\times10^{-4}\ \mathrm{s}$. A similar pattern can be observed: along the $X$-axis pronounced oscillations occur, while along the $Y$-axis the agreement is more accurate.

However, the discrepancies are more pronounced here than in the case of in-plane disturbances. The leading edge of the longitudinal wavefront in the three-dimensional model significantly lags behind the peak observed in the shell model, while the transverse waves exhibit residual oscillations in the vicinity of the initial disturbance.

\begin{figure}[!h]
    \centering
    \begin{minipage}{0.45\linewidth}
        \centering
        \includegraphics[width=\linewidth]{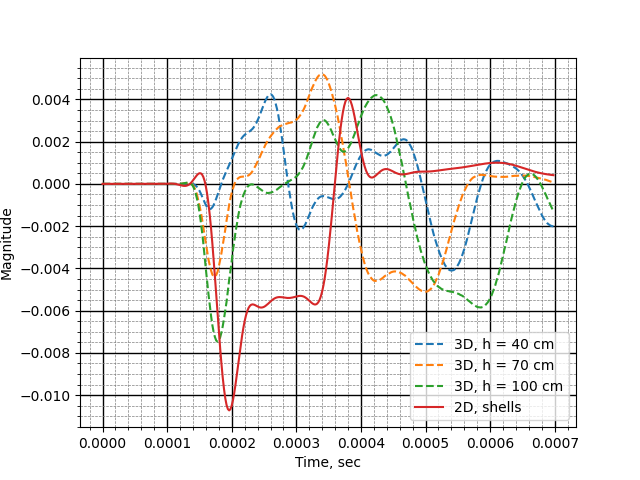}
        \caption{$M_{xx}$ on sensor along X axis}
        \label{fig:sens_Mxx_X}
    \end{minipage}
    \hfill
    \begin{minipage}{0.45\linewidth}
        \centering
        \includegraphics[width=\linewidth]{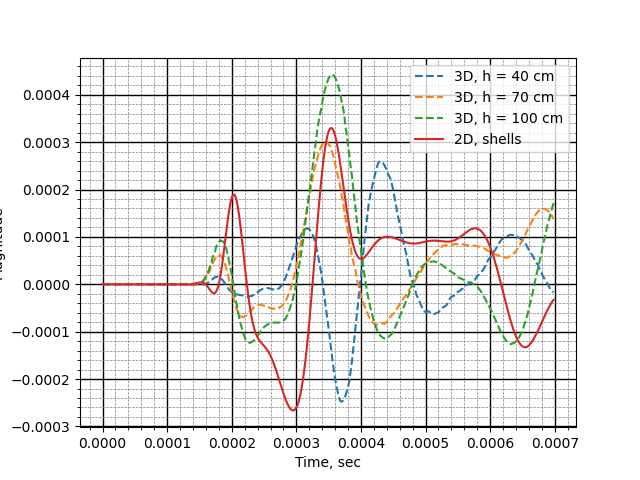}
        \caption{$M_{xx}$ on sensor along Y axis}
        \label{fig:sens_Mxx_Y}
    \end{minipage}
\end{figure}

The sensor simulation results are presented in Figures~\ref{fig:sens_Mxx_X} and \ref{fig:sens_Mxx_Y}. A considerable mismatch is evident: the curves obtained from the three-dimensional model and from the Kirchhoff--Love model do not agree qualitatively.

The NRMSE plot shows that reducing the thickness does not lead to improved accuracy. 
Figure \ref{fig:outplane_sxx_ynorm} shows a significant nonlinearity of wave effects relative to the plane thickness.
It requires further research to establish a critical thickness of the plate for the shell theory to become applicable.
Also, further calculations must be conducted on thinner plates to establish the convergence to the thin shell model.

\begin{figure}[!h]
    \centering
    \includegraphics[width=0.45\linewidth]{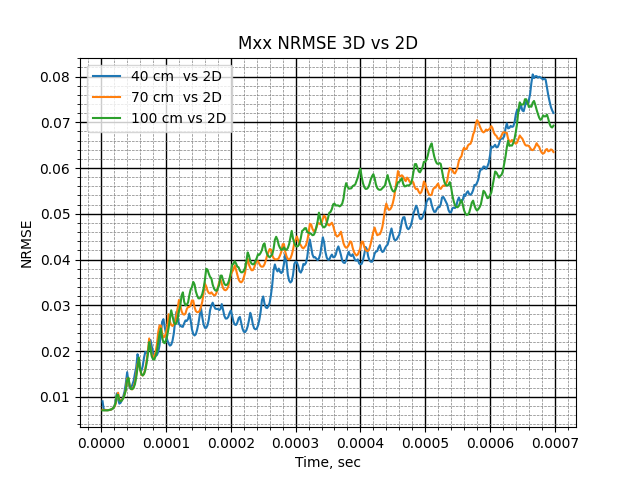}
    \caption{NRMSE metric for $M_{xx}$}
    \label{fig:nrmse_Vx}
\end{figure}

\newpage
\section{Conclusions}

This research on the Kirchhoff-Love theory shows that its dynamic form can be adducted to a hyperbolic set of equations with a substitution that assumes strains and deformations to be small.
It reduces the applicability of this approach to small deformations -- ultrasound, small vibrations and, possibly, brittle damage and fracturing.

The hyperbolic set of equations for Kirchhoff-Love shells separates into two independent parts that describe two independent types of wave processes.
We compared shells results with full three dimensional statements calculated with the grid-characteristic method implementation that was verified in previous research.

The comparison shows a quantitative correspondence for in-plane waves. Reducing the thickness of the plate in full three dimensional statement leads to a better correspondence to the shell model results.

Out-of-plane waves that we consider to be Lamb surface waves show only qualitative correspondence between full three dimensional and shell model results. 
These waves are clearly separate and independent from in-plane waves, and the wave pattern shows the same wave groups. 
However, the qualitative correspondence was not yet achieved and requires further research.

Calculations with the three dimensional model allow to obtain a detailed analysis of wave patterns for both cases and establish specific rules of the applicability of the considered approach to the Kirchhoff-Love shells theory.

\begin{credits}
\subsubsection{\ackname} This study was funded by RSCF (grant number 23-11-00035).

\subsubsection{\discintname}
The authors have no competing interests to declare that are
relevant to the content of this article. 
\end{credits}

\bibliographystyle{hunsrt}
\bibliography{main}

\end{document}